\input amstex
\documentstyle{amsppt}
\magnification=\magstep1
\baselineskip=13pt
\vsize=8.7truein
\CenteredTagsOnSplits
\NoRunningHeads
\def\tr{\operatorname{tr}}
\def\Sym{\operatorname{Sym}}
\def\spa{\operatorname{span}}
\def\En{\operatorname{End}}
\def\FF{\Cal{F}}
\topmatter
\title Estimating Maximum by Moments for Functions 
on Orbits \endtitle
\author Alexander Barvinok \endauthor
\address Department of Mathematics, University of Michigan, Ann Arbor,
MI 48109-1109 \endaddress
\thanks This research was partially supported by NSF Grant DMS 9734138.
 \endthanks
\abstract Let $G$ be a compact group acting in a real vector space
$V$. We obtain a number of inequalities relating the $L^{\infty}$ 
norm of a matrix element of the representation of $G$ with its 
$L^p$ norm for $p<\infty$. We apply our results to obtain  
approximation algorithms to find the maximum absolute value 
of a given multivariate polynomial over the unit sphere (in which 
case $G$ is the orthogonal group) and for the multidimensional 
assignment problem, a hard problem of combinatorial optimization 
(in which case $G$ is the symmetric group). \endabstract
\date January 2002 \enddate
\email barvinok$\@$umich.edu \endemail
\endtopmatter
\document

\head Introduction \endhead

A general problem of optimization has to do with finding the maximum
(minimum)
value of a real valued function $f: X \longrightarrow {\Bbb R}$.
Often, the set $X$ is endowed with a probability measure $\mu$
and the function $f$ possesses a certain degree of symmetry which allows
one to compute the $k$-th moment $\displaystyle \int_X f^k \ d\mu$ 
efficiently at least for small values of $k$. 
Thus one may ask how well the $k$-th moment approximates the maximum value.
In this paper, we 
describe a fairly general situation where some simple and meaningful 
relations
between the maximum and moments can be obtained. 
We provide two illustrations: one, continuous, has to do 
with optimization of multivariate polynomials on the unit sphere with
possible applications to solving systems of real polynomial equations and
the other,
discrete, deals with optimization on the symmetric group, namely, with the 
multidimensional assignment problem, a hard problem of combinatorial 
optimization.    

\subhead (1.1) The general setting \endsubhead
Let $G$ be a compact group with the Haar probability measure $dg$
acting in a finite-dimensional real vector space $V$. 
To avoid dealing with various technical details, we assume that 
the representation $G \longrightarrow GL(V)$ is continuous, where 
the general linear group $GL(V)$ is considered in its standard topology. 

Let us choose a vector $v \in V$ and a linear function 
$\ell : V \longrightarrow {\Bbb R}$. We consider the orbit 
$\bigl\{gv: g \in G \bigr\}$ of $v$ and the resulting function 
$f: G \longrightarrow {\Bbb R}$ defined by
$$f(g)=\ell(gv).$$
In other words, $f$ is a matrix element in the representation of $G$.
We are interested in the relation between the following quantities: 

The $L^{\infty}$ norm of $f$:
$$\|f\|_{\infty}=\max_{g \in G} |f(g)| =\max_{g \in G} \big|
\ell(gv) \big|.$$ 
The $L^{2k}$ norm of $f$ for a positive integer $k$:
$$\|f\|_{2k}=\biggl(\int_{G} f^{2k}(g) \ dg \biggr)^{1 \over 2k}=
\biggl( \int_{G} \ell^{2k}(gv) \ dg \biggr)^{1 \over 2k}.$$
As we remarked earlier, for many examples 
in computational mathematics, the quantity $\|f\|_{\infty}$ is of 
considerable interest and is hard to compute whereas 
$\|f\|_{2k}$ is relatively easy to compute for moderate values of $k$.
First, we relate $\|f\|_{\infty}$ and $\|f\|_2$.

\proclaim{(1.2) Theorem} 
Let $G$ be a compact group acting in a finite-dimensional
real vector space $V$ and let $dg$ be the Haar probability measure 
on $G$. Let us fix a vector $v$ and a linear function 
$\ell: V \longrightarrow {\Bbb R}$
and let us define a real-valued function $f: G \longrightarrow {\Bbb R}$ 
by $f(g)=\ell(gv)$. 
Then 
$$\|f\|_2 \leq \|f\|_{\infty} \leq \sqrt{\dim V}\cdot  \|f\|_2.$$   
\endproclaim

The bounds of Theorem 1.2 are generally sharp, see Remark 2.4. 
To estimate how well $\|f\|_{2k}$ 
approximates $\|f\|_{\infty}$ for a larger $k$,
we invoke a general construction from the representation theory, see 
for example, Lecture 6 of \cite{4}.

\subhead (1.3) Tensor power \endsubhead
For a positive integer $k$, let 
$$V^{\otimes k}=
\underbrace{V \otimes \ldots \otimes V}_{k \text{\ times}}$$
be the $k$-th tensor power of $V$. There is a natural action of $G$ 
in $V^{\otimes k}$, defined on decomposable tensors by
$$g(v_1 \otimes \ldots \otimes v_k)=gv_1 \otimes \ldots \otimes gv_k
\quad \text{for} \quad g \in G.$$ 
There is a natural 
action of the symmetric group $S_k$ permuting the components in the 
tensor product. Thus, for decomposable tensors, we have
$$\sigma(v_1 \otimes \ldots \otimes v_k)=
v_{\sigma^{-1}(1)} \otimes \ldots \otimes v_{\sigma^{-1}(k)} 
\quad \text{for} \quad \sigma \in S_k.$$ 
The action of $S_k$ commutes with the action of 
$G$. Let $\Sym_k(V)$ be the symmetric part of $V^{\otimes k}$ consisting 
of the tensors $x$ such that $\sigma x=x$ for all $\sigma \in S_k$. It 
is known that 
$$\dim \Sym_k(V)={\dim V +k -1 \choose k},$$
since $\Sym_k(V)$ can be thought of as the space of all real homogeneous 
polynomials of degree $k$ in $\dim V$ variables.
Let 
$$v^{\otimes k}=
\underbrace{v \otimes \ldots \otimes v}_{k \text{\ times}}$$
be the $k$-th tensor power of $v$. 
Thus $v^{\otimes k} \in \Sym_k(V)$ and 
$gv^{\otimes k} \in \Sym_k(V)$ for all $g \in G$.
It turns out that how well $\|f\|_{2k}$ approximates $\|f\|_{\infty}$ 
depends on the dimension $D_k$ of the subspace spanned by the orbit 
$\bigl\{g v^{\otimes k}\bigr\}$. This dimension may be different 
for different $v \in V$. Roughly, if $D_k$ is small then $v$ lies in 
a certain algebraic variety constructed from the action of $G$ in $V$ 
and for such $v$ the functions $f$ are ``smoother'' than for those $v$ 
for which $D_k$ is large. 

Thus we obtain the following corollary of Theorem 1.2. 

\proclaim{(1.4) Corollary} Let $G$ be a compact group acting in a 
finite-dimensional
real vector space $V$ and let $dg$ be the Haar probability measure 
on $G$. Let us fix a vector $v$ and a linear function 
$\ell: V \longrightarrow {\Bbb R}$
and let us define a real-valued function $f: G \longrightarrow {\Bbb R}$ 
by $f(g)=\ell(gv)$. For a positive integer $k$, let 
$$D_k=\dim \spa \bigl\{g v^{\otimes k}: g \in G \bigr\}$$ 
be the dimension of the span of the orbit of $v^{\otimes k}$ in 
$V^{\otimes k}$. Then 
$$\|f\|_{2k} \leq \|f\|_{\infty} \leq \bigl(D_k\bigr)^{{1 \over 2k}} 
\cdot \|f\|_{2k}.$$
\endproclaim

Again, generally speaking, the estimates of Corollary 1.4 can 
not be improved, see Remark 2.5.
A straightforward estimate of $D_k \leq \dim \Sym_k(V)$ produces the 
following corollary.

\proclaim{(1.5) Corollary} Let $G$ be a compact group acting in a 
finite-dimensional
real vector space $V$ and let $dg$ be the Haar probability measure 
on $G$. Let us fix a vector $v$ and a linear function 
$\ell: V \longrightarrow {\Bbb R}$
and let us define a real-valued function $f: G \longrightarrow {\Bbb R}$ 
by $f(g)=\ell(gv)$. Let $k$ be a positive integer. Then 
$$\|f\|_{2k} \leq \|f\|_{\infty} \leq 
{\dim V + k -1 \choose k}^{1 \over 2k} \cdot \|f\|_{2k}.$$ 
\endproclaim

There are examples showing that the bounds of Corollary 1.5 are 
``almost tight''. For instance, if $G=SO(n)$ is the orthogonal 
group acting in $V={\Bbb R}^n$, computations of Section 3.3 show 
that the upper bound for $\|f\|_{\infty}$ is tight up to a 
factor of $\sqrt{2}$ (uniformly on $k$ and $n$).
 
As we remarked earlier, in many cases we are able to compute $\|f\|_{2k}$ 
efficiently if $k$ is not very large. Quite often
 (see examples of Sections 3 and 4),
we can compute $\|f\|_{2k}$ in polynomial time for any fixed $k$. 
The following estimate shows the type of bound that we can achieve if 
we fix $k$ in advance.

\proclaim{(1.6) Corollary} For any $\epsilon>0$ there exists a 
$k_0=k_0(\epsilon)=O\bigl(\epsilon^{-2}\bigr)$ 
such that for any positive integer $k>k_0$, for any  
compact group $G$ acting in a real vector space $V$ with $\dim V \geq k$,  
for any linear function $\ell: V \longrightarrow {\Bbb R}$, 
for any $v \in V$ and for the function $f(g)=\ell(gv)$, 
$f: G \longrightarrow {\Bbb R}$, we have 
$$\|f\|_{2k} \leq \|f\|_{\infty} \leq \epsilon \sqrt{\dim V} 
\cdot \|f\|_{2k}.$$  
\endproclaim

The paper is structured as follows. In Section 2, we prove Theorem 1.2 
and Corollaries 1.4--1.6. In Section 3, we apply our results to 
the problem of finding the largest absolute value of a real 
homogeneous multivariate 
polynomial
on the unit sphere, in which case $G=SO(n)$, the orthogonal group.
 In particular, we present a simple polynomial time 
approximation algorithm to compute the largest absolute value on 
the sphere  
of a {\it fewnomial}, that is, a polynomial having only small (fixed) 
number of monomials and discuss possible applications to solving systems
of real fewnomial equations. 
In Section 4, we discuss a hard problem of combinatorial
optimization, that is, the multidimensional assignment problem, in which 
case $G=S_n$. In particular, our results lead to an approximation 
algorithm for finding a bijection between vertex sets of two hypergraphs
$H_1$ and $H_2$, which maximizes the number edges of $H_1$ mapped onto
the edges of $H_2$.

 We use the real model (see \cite{3})
for computational complexity, 
counting the number of arithmetic operations performed by the algorithm. 
Eventually, to compute $\|f\|_{2k}$ from $\|f\|_{2k}^{2k}$ we need to
extract a root of degree $2k$, which we count as a single operation. 

\head 2. Proofs \endhead

In this section, we prove Theorem 1.2 and Corollaries 1.4--1.6. 
We need some standard facts 
from the representation theory (see, for example, \cite{4}).

Let $G$ be a compact group acting in a finite-dimensional real vector 
space $V$. As is known, $V$ possesses a $G$-invariant scalar 
product $\langle \rangle$:
$$\langle u, v \rangle=\langle gu,\ gv \rangle \quad \text{for all} \quad 
u,v \in V \quad \text{and all} \quad g \in G.$$
We introduce the corresponding Euclidean norm:
$$\|x\|=\sqrt{\langle x, x \rangle}.$$
The action (representation) is called irreducible if $V$ contains 
no proper $G$-invariant subspaces. As is known, if $G$ acts in a 
finite-dimensional real vector space $V$, then $V$ can be represented 
as a direct sum of pairwise orthogonal 
(with respect to a given $G$-invariant 
scalar product) invariant subspaces $V_i$ such that the action of $G$ 
in each $V_i$ is irreducible. 

A somewhat ``non-standard'' feature of our construction is that we consider
representations over the real rather than over the complex numbers.
Consequently, we need a substitute for Schur's Lemma. It comes
in the form of the following observation.
   
Suppose that $q: V \longrightarrow {\Bbb R}$ is a $G$-invariant quadratic
form, that is $q(gx)=q(x)$ for all $x \in V$ and all $g \in G$. We 
claim that the eigenspaces of $q$ are $G$-invariant subspaces. 
A possible way  
to see that is to notice that the unit eigenvectors of $q$ are 
precisely the critical points of the restriction 
$q: {\Bbb S} \longrightarrow {\Bbb R}$ where 
${\Bbb S}=\bigl\{x: \|x\|=1 \bigr\}$ is the unit sphere in $V$.

Our first lemma is a real version of the orthogonality relations for 
matrix elements. 
\proclaim{(2.1) Lemma} Let $G$ be a compact group acting in a 
finite-dimensional
real vector space $V$ endowed with a $G$-invariant scalar product 
$\langle \rangle$. Suppose that the representation of $G$ is irreducible 
and let $dg$ be the Haar probability measure on $G$.
Then
$$\int_G \langle x,\ gv \rangle^2 \ dg =
{\|v\|^2 \cdot \|x\|^2 \over \dim V} \quad 
\text{for all} \quad x, v \in V.$$
\endproclaim 
\demo{Proof}
Let us choose a vector $v \in V$ and let us define a quadratic form 
$q: V \longrightarrow {\Bbb R}$ by 
$$q(x)=\int_G \langle x,\ gv \rangle^2 \ dg.$$
Clearly, $q(x)$ is $G$-invariant: $q(gx)=q(x)$ for all $x \in V$ 
and all $g \in G$.
Let $\lambda$ be the largest eigenvalue of $q$ and let $W$ 
be the corresponding 
eigenspace. Then $W$ is an invariant subspace of $V$ and hence $W=V$.
Thus $q(x)=\lambda\|x\|^2$ for some $\lambda \geq 0$. 

To find $\lambda$, let us compute the trace of $q$.
On one hand, we have $\tr q=\lambda \dim V$. 

Let $q_g(x)=\langle x,\ gv \rangle^2$. Then $q_g$ is a quadratic form 
of rank 1 with the non-zero eigenvalue $\|gv\|^2=\|v\|^2$
which corresponds to an eigenvector $x=gv$. Hence $\tr q_g=\|v\|^2$ 
and since $q(x)$ is the average of $q_g$, we have $\tr q= \|v\|^2$.
Thus $\lambda=\|v\|^2/\dim V$. Hence 
$$q(x)={\|v\|^2 \cdot \|x\|^2 \over \dim V}$$ 
and the proof follows.
{\hfill \hfill \hfill} \qed
\enddemo
Now we use that every representation is a sum of irreducible 
representations.
\proclaim{(2.2) Lemma} 
Let $G$ be a compact group acting in a finite-dimensional
real vector space $V$ endowed with a $G$-invariant scalar product 
$\langle \rangle$. Let $dg$ be the Haar probability measure on $G$.
Let us fix a vector $v \in V$. Then there exists a decomposition 
$V=V_1 \oplus \ldots \oplus V_k$ of $V$ into the direct sum of non-zero 
pairwise orthogonal invariant subspaces such that for every $x \in V$
we have 
$$\int_G \langle x,\ gv \rangle^2 \ dg =
\sum_{i=1}^k {\|x_i\|^2 \cdot \|v_i\|^2 \over \dim V_i},$$
where $x_i$ and $v_i$ are the orthogonal projections onto $V_i$ of $x$ and
$v$ respectively. 
\endproclaim
\demo{Proof}
Let us define a quadratic form $q: V \longrightarrow {\Bbb R}$ by 
$$q(x)=\int_{G} \langle x,\ gv \rangle^2 \ dg.$$
Then $q$ is $G$-invariant, $q(gx)=q(x)$ for all $g \in G$ and all 
$x \in V$. Thus the eigenspaces of $q$ are $G$-invariant subspaces of $V$.
 Let us 
write every eigenspace as a direct sum of pairwise orthogonal 
invariant subspaces $V_i$ such that the action of $G$ in each $V_i$ 
is irreducible. Thus we obtain the decomposition 
$V=V_1 \oplus \ldots \oplus V_k$ and we have  
$$q(x)=\sum_{i=1}^k \lambda_i \|x_i\|^2,$$
where $x_i$ is the orthogonal projection of $x$ onto $V_i$ and $\lambda_i$
are non-negative numbers.
To find $\lambda_i$, let us choose a non-zero $x \in V_i$.  
Then $\langle x, g v \rangle=\langle x, g v_i \rangle$ and by Lemma 2.1,
we get 
$$q(x)=\lambda_i \|x\|^2 ={\|v_i\|^2 \cdot \|x\|^2 \over \dim V_i},$$
from which 
$$\lambda_i={\|v_i\|^2 \over \dim V_i}.$$ 
The proof now follows.
{\hfill \hfill \hfill} \qed
\enddemo
\remark{(2.3) Remark} A decomposition $V=V_1 \oplus \ldots \oplus V_k$ of 
a representation into the direct sum of pairwise orthogonal 
irreducible components is {\it not}
unique as long as some irreducible representation appear with a 
multiplicity greater than 1 (which means that the representations of $G$ 
in some subspaces $V_i$ are isomorphic). One can construct some simple 
examples showing that the decomposition of Lemma 2.2 indeed depends 
on $v$. 
\endremark
Now we are ready to prove Theorem 1.2.

\demo{Proof of Theorem 1.2} The inequality 
$$\|f\|_2 \leq \|f\|_{\infty}$$ 
is quite standard. Let us prove that
$$\|f\|_{\infty} \leq \sqrt{\dim V} \cdot \|f\|_2.$$ 
Let $e$ be the identity in $G$.
We note that it suffices to prove that 
$$|f(e)|=|\ell(v)| \leq \sqrt{\dim V} \cdot \|f\|_2,$$
because the inequality for $f(g)=\ell(gv)$ would follow then 
by choosing a new vector $v$: 
$$\text{new\ }v:=g \bigl(\text{\ old \ }v\bigr).$$ 

Let us introduce a $G$-invariant scalar product 
$\langle \rangle$ in $V$ so that $\ell(x)=\langle c, x \rangle$ for some 
$c \in V$ and all $x \in V$. Applying Lemma 2.2, we obtain a decomposition
$V=V_1 \oplus \ldots \oplus V_k$ into the direct sum of pairwise 
orthogonal invariant subspaces such that 
$$\|f\|_2^2=\int_{G} \langle c,\ gv \rangle^2 \ dg =
\sum_{i=1}^k {\|c_i\|^2 \cdot \|v_i\|^2 \over \dim V_i},$$
where $c_i$ and $v_i$ are the orthogonal projections onto $V_i$ of $c$ and 
$v$ 
respectively.
We have 
$$f(e)=\langle c, v \rangle=\sum_{i=1}^k \langle c_i, v_i \rangle$$
and hence 
$$ |f(e)| \leq \sum_{i=1}^k |\langle c_i, v_i \rangle| \leq 
\sum_{i=1}^k \|c_i\| \cdot \|v_i\|.$$
Let 
$$\alpha_i={\|c_i\| \cdot \|v_i\| \over \sqrt{\dim V_i}}.$$
Then 
$$|f(e)|^2  \leq \biggl(\sum_{i=1}^k \alpha_i \sqrt{\dim V_i}\biggr)^2 
\leq
\Bigl( \sum_{i=1}^k \alpha_i^2 \Bigr) \Bigl( \sum_{i=1}^k \dim V_i \Bigr)=
\bigl(\dim V \bigr) \cdot \|f\|_2^2.$$
and the proof follows.  
{\hfill \hfill \hfill} \qed
\enddemo
\remark{(2.4) Remark} Analyzing the proof of Theorem 1.2, it is not hard 
to find out when the bound 
$\|f\|_{\infty} \leq \sqrt{\dim V} \cdot \|f\|_2$ is sharp. In particular,
the bound is sharp for the class of linear functions on the 
orbit of $v$ as long as the orbit of $v$ spans $V$.  
Here are some natural cases when the bound is attained. 

Suppose, for example, that we have an absolutely irreducible representation
$\rho$ of $G$ in a real vector space $W$ 
(that is, the representation remains irreducible after 
complexification). Thus, for every $g \in G$, 
$\rho(g)$ is an operator 
in $W$. We interpret $\rho(g)$ as a point in the space 
$V=\En(W)$ of all linear transformations $W \longrightarrow W$. Let 
$\chi(g)=\tr(g)$ be the character of the representation. We  
think of $\chi(g)$ as of a linear function on the orbit of the identity 
operator $I \in \En(W)$ under the action $g(x)=\rho(g)x$ for all 
$x \in \En(W)$. Then $\|\chi\|_{\infty}=\dim W=\sqrt{\dim V}$. 
The orthogonality relations for the characters  (see, for example, 
Chapter 2 of \cite{4}) state that 
$\|\chi\|_2=1$ and hence the bound of 
Theorem 1.2 holds with equality.

As another example, let us consider a finite group $G$ 
of cardinality $|G|$ and an arbitrary
function $f: G \longrightarrow {\Bbb R}$. Of course, in this case,
the inequality $ \|f\|_{\infty} \leq \sqrt{|G|} \cdot \|f\|_2$ 
is the best we can hope for (take $f$ to be the 
delta-function of an element of $G$). The function $f$ can be thought of
as 
a linear function on the orbit of a point in the regular representation 
of $G$. The space $V$ in this case is the vector space of all 
linear functions $f: G \longrightarrow {\Bbb R}$ where $G$ acts by 
shifts: $gf(x)=f(g^{-1}x)$. 
Let $v \in V$ be the delta-function at the identity: 
$v(e)=1$ where $e$ is the identity in $G$ and $v(g)=0$ for all $g \ne e$.
Then $f$ is a linear function on the orbit of $v$ and $\dim V=|G|$.
\endremark
\bigskip
To prove Corollary 1.4, we use the construction of the tensor power 
(see Section 1.3). 
\demo{Proof of Corollary 1.4}
Let us define a function $h: G \longrightarrow {\Bbb R}$ by 
$$h(g)=f^k(g)=\ell^{\otimes k} (v^{\otimes k}).$$
Thus $h$ is a linear function on the orbit of $v^{\otimes k}$.
Let 
$$W=\spa\bigl\{g v^{\otimes k} : \ g \in G \bigr\}$$
be the span of the orbit of $v^{\otimes k}$. Hence $\dim W=D_k$.
Applying Theorem 1.2 to the linear function $h$ on the orbit 
of $v^{\otimes k}$ in $W$, we get 
$$\|h\|_2 \leq \|h\|_{\infty} \leq \sqrt{D_k} \cdot \|h\|_2.$$
Now we note that
$\|h\|_{\infty}=\|f\|_{\infty}^k$ and that
$\|h\|_2=\|f\|_{2k}^k$.
{\hfill \hfill \hfill} \qed
\enddemo

\remark{(2.5) Remark} The bound 
$\|f\|_{\infty} \leq \bigl(D_k\bigr)^{1 \over 2k} \|f\|_{2k}$ is rarely 
sharp. One example when it is sharp is provided by a generic 
orbit in the regular representation of a finite group $G$,
see Remark 2.4. In Section 3.3, we present a series of 
examples of matrix functions 
for $G=SO(n)$ for which the estimate is sharp up to a constant factor
uniformly on $k$ (and uniformly on $n$).   
\endremark
\bigskip
Corollary 1.5 follows by a general estimate of $D_k$.
\demo{Proof of Corollary 1.5}
We apply Corollary 1.4. The orbit $\bigl\{gv^{\otimes k}\bigr\}$ lies 
in the 
symmetric part $\Sym_k(V)$ of the tensor product $V^{\otimes k}$ and 
hence 
$$D_k=\dim \spa \bigl\{g v^{\otimes k} \bigr\} \leq \dim \Sym_k(V) 
={\dim V + k -1 \choose k}.$$
{\hfill \hfill \hfill} \qed
\enddemo

\remark{(2.6) Remark} As follows from Section 3.3, the 
upper bound for $\|f\|_{\infty}$ is sharp up to 
a constant factor uniformly on $k$ if $G=SO(n)$, $V={\Bbb R}^n$ and
$G$ acts in $V$ by its defining representation. 

We describe below classes of functions $f: G \longrightarrow {\Bbb R}$
for which some sharp estimates can be obtained.
Let us fix a representation $\rho$ of $G$ in a real 
vector space $V$ and let $\FF_{\rho}$ be the vector space spanned by the 
matrix elements of $\rho$. Then for some constant 
$C(\rho, k)$ and for all $f \in \FF_{\rho}$ we have 
$\|f\|_{\infty} \leq C(\rho, k) \|f\|_{2k}$.
Let us assume that $\rho$ is absolutely irreducible
(cf. Remark 2.4). In principle, the best
possible value of $C(\rho, k)$ can be computed from the representation 
theory of $G$ as follows. Let us choose an $f \in \FF_{\rho}$. 
Shifting $f$,
if necessary, we may assume that the maximum absolute value of $f$ is
attained at the identity $e$ of $G$. Let us define 
$h: G \longrightarrow {\Bbb R}$ by 
$$h(x)=\int_{G} f(g^{-1} x g) \ dg \quad \text{for all} \quad x \in G.$$
Then $\|h\|_{\infty}=\|f\|_{\infty}$ and $\|h\|_{2k} \leq \|f\|_{2k}$
for all positive integers $k$. Thus the largest ratio 
$\|f\|_{\infty}/\|f\|_{2k}$ for $f \in \FF_{\rho}$ is attained when 
$f$ satisfies $f(g^{-1} x g)=f(x)$ for all $g \in G$ and all $x \in G$
from which it follows that $f$ is a multiple of the character 
$\chi(g)=\tr \rho(g)$,
see Remark 2.4. We observe that $\|\chi\|_{\infty}=\dim V$. 
Moreover, the orthogonality relations (see Lecture 6 of \cite{4}) imply that 
$\|\chi\|_{2k}^{2k}$ is the sum of squares of multiplicities 
of the irreducible components of the tensor power $\rho^{\otimes k}$.
Summarizing, we conclude that
 in order to be able to compute the best possible constant 
$C(\rho, k)$ such that $\|f\|_{\infty} \leq C(\rho, k) \|f\|_{2k}$ for 
any linear combination $f$ of matrix elements of $\rho$, 
it suffices to know how the tensor power 
$\rho^{\otimes k}$ decomposes into the sum of absolutely irreducible 
representations.   
\endremark
\bigskip
Finally, Corollary 1.6 follows by an estimate of the binomial coefficient.
\demo{Proof of Corollary 1.6}
Let us choose a $k_0$ such that $(k!)^{1/k}>2\epsilon^{-2}$
for all $k>k_0$. 
By Stirling's formula, 
we can choose $k_0=O\bigl(\epsilon^{-2}\bigr)$.
Then 
$$\split {\dim V + k -1 \choose k}^{1 \over 2k} &=
\biggl({\dim V \cdot (\dim V +1) \cdots (\dim V +k-1) \over k!}
\biggr)^{1\over 2k} \\ &\leq 
\biggl({2^k \dim^k V \over k!}\biggr)^{1 \over 2k}=
2^{1/2} (k!)^{-1/2k} \cdot \sqrt{\dim V} \\
& \leq \epsilon \sqrt{\dim V}. \endsplit$$
The proof follows by Corollary 1.5. 
{\hfill \hfill \hfill} \qed
\enddemo

\head 3. Applications to Polynomials \endhead

In this section, we apply our results to approximate the maximum 
absolute value of a homogeneous multivariate polynomial on the 
unit sphere. 

Let $p$ be a homogeneous polynomial of degree $d$ in $n$ real variables 
$\xi_1, \ldots, \xi_n$. Thus we can write
$$p(x)=\sum_{1 \leq i_1, \ldots, i_d \leq n} \gamma_{i_1 \ldots i_d}
\xi_{i_1} \cdots \xi_{i_d} \quad \text{for} \quad 
x=(\xi_1, \ldots, \xi_n)$$
where $\gamma_{i_1 \ldots i_d}$ are some real numbers.  

Let ${\Bbb R}^n$ be the $n$-dimensional Euclidean space and 
let $x=(\xi_1, \ldots, \xi_n) \in {\Bbb R}^n$ be a point. Then 
$V=\bigl({\Bbb R}^n\bigr)^{\otimes d}$ can be identified with the space 
${\Bbb R}^{n^d}$.
The coordinates of a typical point (tensor) $X \in V$ are 
$$\Bigl(X_{i_1 \ldots i_d}: \quad 1 \leq i_1, \ldots, i_d \leq n \Bigr)$$
and the scalar product in $V$ is defined by 
$$\big\langle X,\ Y \big\rangle=
\sum_{1 \leq i_1, \ldots, i_d \leq n} 
X_{i_1 \ldots i_d} Y_{i_1 \ldots i_d}.$$ 
For $x=(\xi_1, \ldots, \xi_n) \in {\Bbb R}^n$, 
the coordinates of $x^{\otimes d}$
are 
$$\Bigl( \xi_{i_1} \cdots \xi_{i_d} \quad \text{for} \quad 
1 \leq i_1, \ldots, i_d \leq n \Bigr).$$
Therefore, we can write 
$$p(x)=\big\langle c,\ x^{\otimes d} \big\rangle \quad 
\text{where} \quad c=\Bigl(\gamma_{i_1 \ldots i_d} \Bigr).$$ 
Let $G=SO(n)$ be the group of orientation preserving orthogonal 
transformations of ${\Bbb R}^n$. Then $G$ acts in $V$ by 
the $d$-th tensor power of its defining representation in ${\Bbb R}^n$.
Let us choose 
$w=(1,0, \ldots, 0) \in {\Bbb R}^n$.
Then, for any $g \in G$, we have 
$$\big\langle c,\ g w^{\otimes d} \big\rangle =p(gw)$$ 
and the orbit 
$\{gw: g \in G \}$ is the unit sphere ${\Bbb S}^{n-1} \subset {\Bbb R}^n$.
Thus the values of $p(x)$, as $x$ ranges over the unit sphere in 
${\Bbb R}^n$, are the values of the linear function 
$$\ell(g w^{\otimes d})=\langle c,\ g w^{\otimes d} \rangle=
\langle c,\ x^{\otimes d} \rangle$$
as $g$ ranges over the orthogonal group $SO(n)$.
 
Moreover, the push-forward of the Haar probability measure $dg$ on $G$ 
is the probability measure $dx$ on ${\Bbb S}^{n-1}$.
Thus we connect the values of a polynomial on the unit sphere with 
the values of a linear function on the orbit of the group $G=SO(n)$.

\proclaim{(3.1) Corollary} Let $p$ be a homogeneous polynomial of degree 
$d$ 
in $n$ real variables, let ${\Bbb S}^{n-1}$ be the unit sphere in 
${\Bbb R}^n$ and let $dx$ be the rotation invariant probability measure 
on ${\Bbb S}^{n-1}$. For a positive integer $k$, let us define the 
$L^{2k}$ norm of $p$ by
$$\|p\|_{2k}=
\biggl( \int_{{\Bbb S}^{n-1}} p^{2k}(x) \ dx \biggr)^{1\over 2k}$$ 
and the $L^{\infty}$ norm by 
$$\|p\|_{\infty}=\max_{x \in {\Bbb S}^{n-1}} |p(x)|.$$
Then 
$$\|p\|_{2k} \leq \|p\|_{\infty} \leq {kd+n-1 \choose kd}^{1 \over 2k}
\|p\|_{2k}.$$ 
\endproclaim
\demo{Proof} We apply Corollary 1.5. Let $w=(1,0, \ldots, 0)$
be as above. Then, for $v=w^{\otimes d}$, we can write 
$p(gw)=\ell(gv)$ for some linear functional 
$\ell: V \longrightarrow {\Bbb R}$ and all 
$g \in G$. The dimension $D_k$ of the 
span of the orbit $\bigl\{gv^{\otimes k}=gw^{\otimes kd}: \ g \in G \bigr\}$ 
is that of the space of 
homogeneous polynomials of degree $kd$ in $n$ variables. 
Hence 
$$D_k={kd+n-1 \choose kd}.$$
We have 
$$\int_{G} \ell^{2k}(gv) \ dg =\int_{{\Bbb S}^{n-1}} p^{2k}(x) \ dx.$$ 
The proof now follows.   
{\hfill \hfill \hfill} \qed
\enddemo 
One way to integrate polynomials over the unite sphere is to 
take the sum of the integrals of the monomials. The following result is 
certainly known, but for the sake of completeness, we sketch its 
proof here.

\proclaim{(3.2) Lemma} Let $p(x)=\xi_1^{\alpha_1} \cdots 
\xi_n^{\alpha_n}$ be 
a monomial. If at least one of $\alpha_i$'s is odd then 
$$\int_{{\Bbb S}^{n-1}} p(x) \ dx=0.$$
If $\alpha_i=2\beta_i$, where $\beta_i$ are non-negative integers 
for $i=1, \ldots, n$, then 
$$\int_{{\Bbb S}^{n-1}} p(x) \ dx = 
{\Gamma(n/2) \prod_{i=1}^n \Gamma(\beta_i + 1/2) \over \pi^{n/2} 
\Gamma\bigl(\beta_1 + \ldots + \beta_n + n/2\bigr)},$$
where $dx$ is the Haar probability measure on ${\Bbb S}^{n-1}$.
\endproclaim
\demo{Sketch of Proof} If $\alpha_i$ is odd then 
$$p(\xi_1, \ldots, \xi_{i-1}, -\xi_i, \xi_{i+1}, \ldots,  \xi_n)=
-p(\xi_1, \ldots, \xi_{i-1}, \xi_i, \xi_{i+1}, \ldots,  \xi_n)$$ 
and hence the average value of $p$ over the unit sphere is 0.

Assuming that $\alpha_i=2 \beta_i$ for $i=1, \ldots, n$, we get 
$$\int_{{\Bbb R}^n} p(x) e^{-\|x\|^2} \ d\mu =
\prod_{i=1}^n \int_{\Bbb R} \xi^{2\beta_i} e^{-\xi^2} \ d\xi=
\prod_{i=1}^n \Gamma\bigl(\beta_1 +1/2\bigr),$$
where $\mu$ is the standard Lebesgue measure in ${\Bbb R}^n$.

On the other hand, passing to the polar coordinates and using that
$p$ is homogeneous of degree $d=2(\beta_1, + \ldots  +\beta_n)$, we 
get
$$\int_{{\Bbb R}^n} p(x) e^{-\|x\|^2} \ d\mu=
|{\Bbb S}^{n-1}| \cdot  
\Bigl(\int_{{\Bbb S}^{n-1}} p(x) \ dx\Bigr) 
\cdot \Bigl( \int_0^{+\infty} r^{d+n-1} e^{-r^2} \ dr \Bigr),$$
where 
$|{\Bbb S}^{n-1}|=2 \pi^{n/2}/\Gamma(n/2)$ is the Euclidean volume 
of the unit sphere in ${\Bbb R}^n$.
The proof now follows.   
{\hfill \hfill \hfill} \qed
\enddemo
The estimates of Corollary 3.1 are probably not optimal (apart from 
the case of $k=1$), but the following simple example shows that in some
sense, they are close to being optimal.

\subhead (3.3) Powers of linear functions \endsubhead
Let $p$ be the power of a linear function, for example,
 $p(x)=\xi_1^d$. Then 
$\|p\|_{\infty}=1$ and, by Lemma 3.2, 
$$\|p\|_{2k}=
\biggl( {\Gamma(n/2) \Gamma(kd+1/2) \over \sqrt{\pi} \Gamma(kd +n/2)}
\biggr)^{1 \over 2k}.$$ 
Then Corollary 3.1 gives us the estimate 
$$\split \|p\|_{\infty} &\leq \biggl({\Gamma(n/2) \Gamma(kd+1/2) 
\Gamma(kd +n) \over \sqrt{\pi} \Gamma(kd+n/2) \Gamma(n) \Gamma(kd+1)}
\biggr)^{1 \over 2k} 
 \leq \biggl({\Gamma(n/2) \Gamma(kd+n) \over \Gamma(kd + n/2) \Gamma(n)}
\biggr)^{1 \over 2k} \\
&=\biggl({n(n+1) \cdots (kd+n-1) \over (n/2) (n/2+1) \cdots (kd+n/2-1)}
\biggr)^{1 \over 2k} \leq 2^{d \over 2}. \endsplit$$  
Hence, among all homogeneous polynomials of a given degree $d$,
powers of linear functions give the largest ratio
$\|f\|_{\infty}/\|f\|_{2k}$ up to a constant factor depending on 
the degree of $f$ and independent of the number of variables 
$n$ and the value of $k$. This may serve as an indication that the 
bound of Corollary 1.5 are not too bad, cf. Remarks 2.5 and 2.6. 
G. Blekherman \cite{2}
 pointed out to the author that the powers, in general, 
do not provide exactly the largest ratio $\|f\|_{\infty}/\|f\|_{2k}$
among all polynomials of a given degree $d$.
Such ``extremal'' polynomials $f$ were computed by G. Blekherman 
when some of the parameters $n$, $d$ and $k$ are small.
\bigskip
Suppose we want to approximate $\|p\|_{\infty}$ by $\|p\|_{2k}$ for 
a sufficiently large $k$. Let us see what trade-off between 
between the computational complexity and accuracy can we achieve.

\subhead (3.4) Low degree polynomials \endsubhead
Let us fix the degree $d$ and allow the number of variables to vary.
Suppose 
that we are given a 
homogeneous polynomial $p$ of degree $d$ and that we want to estimate 
$\|p\|_{\infty}$. This problem is provably computationally hard 
already for $d=4$ (one can infer it from results 
of Part 1 of \cite{3}) and is suspected to be hard for $d=3$.

Let $m$ be the number of monomials in $p$, so $m=O(n^d)$.
We observe that for any fixed $k$, the direct computation of $p^{2k}(x)$ 
and computing $\|p\|_{2k}$ via Lemma 3.2 has $O(m^{2k})$ complexity.
One the other hand, using Corollary 3.1, we get that 
$$\|p\|_{\infty} \leq C(k)n^{d/2} \|p\|_{2k} \quad \text{where} \quad 
C(k)=O\bigl(k^{-1/2}\bigr).$$
In other words, for any fixed $\epsilon>0$ there is a polynomial time 
algorithm estimating $\|p\|_{\infty}$ within 
a factor of $\epsilon n^{d/2}$. 
If we want a better estimate, we have to take a larger $k$.
Thus, for any constant $C>1$, from Corollary 3.1 (cf. also 
Corollary 1.6), we get that 
$$\|p\|_{\infty} \leq C \|p\|_{2k} \quad \text{for some} \quad 
k=O(n).$$
Since $p^{2k}(x)$ contains at most ${2kd +n -1 \choose 2kd}$ monomials,
we can compute $\|p\|_{2k}$ by Lemma 3.2 in $2^{O(n)}$ time.
Summarizing, for any $C>1$ there exists a $\gamma>0$ such that we can 
approximate $\|p\|_{\infty}$ within a factor $C$ in
$2^{\gamma n}$ time.

\subhead (3.5) Fewnomials \endsubhead
Suppose that we do not fix the degree $d$ of $p$ but 
fix instead the number $m$ of monomials in $p$.
Thus we can write 
$$p(x)=\sum_{i=1}^m p_i(x),$$
where 
$$p_i(x)=\gamma_i \xi_1^{\alpha_{i1}} \cdots \xi_n^{\alpha_{in}}$$
are monomials.
For a positive integer $k$, by the multinomial expansion, we get 
$$p^{2k}=\sum \Sb r_1, \ldots, r_m \geq 0 \\
r_1 + \ldots + r_m = 2k \endSb 
{(2k)! \over r_1! \cdots r_m!} p_1^{r_1} \cdots p_m^{r_m}.$$
 Thus 
$p^{2k}$ contains at most ${m+2k-1 \choose m-1}$ monomials, which is 
a polynomial in $k$ when $m$ is fixed.
Using Lemma 3.2, we compute $\|p\|_{2k}$ in $O\bigl(dn(2k)^m\bigr)$ time.
Given an $\epsilon>0$, let us choose an integer 
$k=O(\epsilon^{-1} n^2 \ln d)$ such 
that 
$${n-1 \over 2k} \ln \bigl(kd+1\bigr) < \ln (1+\epsilon).$$
Using Corollary 3.1 and a simple estimate 
$${kd +n -1 \choose n-1}=
{(kd+1)(kd+2) \cdots (kd+n-1) \over 1 \cdot 2 \cdots (n-1)}  
\leq (kd+1)^{n-1}$$
we conclude that 
$$\|p\|_{2k} \leq \|p\|_{\infty} \leq (1+\epsilon)\cdot \|p\|_{2k}.$$
Hence as long as the number of monomials is fixed, we get a 
polynomial time approximation algorithm, which, for any 
given $\epsilon>0$ computes the maximum 
absolute value of a given polynomial (``fewnomial'') over the 
unit sphere within a relative error of $\epsilon$, 
in time polynomial in $\epsilon^{-1}$, the number of variables $n$ and 
the degree $d$ of the polynomial. In fact, the only place where we have 
to use polynomially many in $d$ arithmetic operations is when we compute 
gamma-functions (factorials) in Lemma 3.2. Apart from that, the
running time of the algorithm is 
polynomial in $\ln d$.
\bigskip
Computing or approximating the maximum absolute value of a polynomial 
on the unit sphere can be used for testing whether a given system 
of real polynomial equations has a real solution, a difficult and
important problem, see for example, \cite{3} and \cite{7}. Suppose that 
$p_i$: $i=1, \ldots, s$ 
are given homogeneous polynomials of degree $d$ in $n$ variables
$x=(\xi_1, \ldots, \xi_n)$ and that we would like 
to test whether the system 
$$p_i(x)=0 \quad \text{for} \quad i=1, \ldots, s$$ 
has a real solution $x \ne 0$. 
Let 
$$q=\sum_{i=1}^s p_i^2(x).$$
Thus we want to test whether 
$$\min_{x \in {\Bbb S}^{n-1}} q(x)=0.$$
Let us choose a 
$$\gamma > \max_{x \in {\Bbb S}^{n-1}} q(x)$$
and let 
$$p=\gamma \|x\|^{2d} -q.$$
Thus the problem reduces to checking whether
$$\max_{x \in {\Bbb S}^{n-1}} |p(x)|=\gamma.$$
If the polynomials $p_i$ of the original system did not have
too many monomials,
we can try to approximate $\|p\|_{\infty}$ by $\|p\|_{2k}$ for 
a reasonably large $k$, cf. Section 3.5. 
Similarly, to choose an appropriate $\gamma$, we can compute 
$\|q\|_{2k}$ for a sufficiently large $k$.
The number of monomials in the 
system is relevant to the ``topological 
complexity'' of the set of real solutions \cite{5}, so it should 
not be surprising that it is also relevant to the computational 
complexity of the decision problem.  
In particular, this approach 
may be useful for detecting ``badly unsolvable'' systems 
(systems for which the value of $\|p\|_{\infty}$ is substantially smaller
than $\gamma$) of fewnomial equations.

\head 4. Applications to Combinatorial Optimization \endhead

Let us fix a number $d$ and let 
$V=\bigl({\Bbb R}^n\bigr)^{\otimes d}={\Bbb R}^{n^d}$ be the vector 
space of $d$-dimensional arrays (tensors) 
$$X=\Bigl(x_{i_1 \ldots i_d}: 
\quad 1 \leq i_1, \ldots, i_d \leq n \Bigr).$$
To simplify the notation somewhat, we denote the coordinates of 
$X$ by $x_I$, where $I=(i_1, \ldots, i_d)$.

We introduce the scalar product by
$$\big\langle X, Y \big\rangle=\sum_{I} x_I y_I \quad 
\text{for} \quad I=(1 \leq i_1, \ldots, i_d \leq n).$$
Let $G=S_n$ be the symmetric group of all permutations $g$ of the 
set $\{1, \ldots, n\}$. We introduce the action of $S_n$ on $V$ 
by the $d$-th tensor power of the natural action of $S_n$ in 
${\Bbb R}^n$:
$$Y=g X \quad \text{provided} \quad 
x_I=y_{g I} \quad \text{where} \quad 
g\bigl(i_1, \ldots, i_d\bigr)=
\bigl(g(i_1), \ldots, g(i_d)\bigr).$$
Let us choose two tensors $A, B \in  V$ and let 
$$f(g)=\big\langle B,\ g A \big\rangle=
\sum_{1 \leq i_1, \ldots, i_d \leq n}
a_{i_1 \ldots i_d} b_{g(i_1) \ldots g(i_d)}, \qquad 
f: S_n \longrightarrow {\Bbb R}$$ be the corresponding matrix 
element.

The problem of maximizing (minimizing) $f$ is one of the 
most general problems of combinatorial optimization, known as the 
$d$-{\it dimensional assignment problem} (see, for example, \cite{6}).
 It is straightforward for 
$d=1$ but already quite difficult for $d=2$ (see \cite{1}).

\example{(4.1) Example: hypergraphs} Recall that a 
$d$-{\it hypergraph} $H$ on the set 
\newline $\{1, \ldots, n\}$ is a set 
of subsets $E \subset \{1, \ldots, n\}$,
called {\it edges} of $H$, such that $|E| \leq d$ for 
the cardinality $|E|$ of every edge $E$ of $H$. 
A hypergraph is called {\it uniform} provided
$|E|=d$ for every edge $E$ of $H$.
Let $H_1$ and $H_2$ be uniform 
$d$-hypergraphs with the set of vertices $\{1, \ldots, n\}$. Let us 
define the {\it adjacency tensor} $A=(a_{i_1 \ldots i_d})$ of $H_1$
by $$a_{i_1 \ldots i_d}
=\cases 1 &\text{if} \quad \{i_1, \ldots, i_d\} \quad \text{is an edge of} 
\quad H_1 \\ 0 &\text{otherwise.} \endcases$$
Let us define $B=(b_{i_1 \ldots i_d})$ by: 
$$b_{i_1 \ldots i_d}=
\cases {1 \over d!} 
&\text{if} \quad \{i_1, \ldots, i_d\} \quad \text{is an edge of}
\quad H_2 \\ 0 &\text{otherwise.} \endcases$$
A permutation $g$ of the set $\{1, \ldots, n\}$ is interpreted 
as a bijection between the vertices of $H_2$ and the vertices of $H_1$ 
and the value of
$$f(g)=\big\langle B,\ gA \big\rangle$$ is the 
number of edges of $H_2$ mapped onto the edges of $H_1$. 
The value of $\|f\|_{\infty}$ is the maximum number of edges of $H_1$ and 
$H_2$ that can be matched by a bijection of the vertices of $H_1$ and
$H_2$. If $H_1$ and $H_2$ are not uniform, we can modify 
$B$ by letting 
$$b_{i_1, \ldots, i_d}={k_1 ! \cdots k_r! \over d!}$$
provided $\{i_1, \ldots, i_d\}$ is an edge of $H_2$ and the 
multiplicities of the elements in the multiset $\{\{i_1, \ldots, i_d\}\}$ 
are $k_1, \ldots, k_r$, so that $k_1 + \ldots +k_r=d$.
Then again the value of $\|f\|_{\infty}$ is equal to the maximum 
number of edges of $H_1$ and $H_2$ that can be matched by a bijection 
of the vertex sets. 

 One can extend this construction to oriented 
hypergraphs whose edges are {\it ordered} subsets of 
$\{1, \ldots, n\}$.
By introducing weights on the edges of $H_1$ and $H_2$ we can
introduce ``prices'' for matching (or mismatching) particular edges.
\endexample 

Applying Corollary 1.5, we get the inequality
$$\|f\|_{2k} \leq \|f\|_{\infty} \leq 
{n^d + k -1 \choose k}^{1 \over 2k} \|f\|_{2k}$$
for the function $f$ of a general $d$-dimensional assignment problem.

In various special cases, the bound can be somewhat improved by
using Corollary 1.4. 
For example, if the coordinates of $A$ (or $B$) are 0's and 1's,
one can prove that   
$$\|f\|_{2k} \leq \|f\|_{\infty} 
\leq D^{1 \over 2k}(n,d,k) \cdot \|f\|_{2k} \quad 
\text{where} \quad D(n,d,k)=\sum_{j=1}^k {n^d \choose j}.$$
 
We claim that for small (fixed) 
values of $k$ the value of $\|f\|_{2k}$ can be computed 
relatively easily (in polynomial time).
First, we observe that computation of $\|f\|_{2k}$ reduces 
to computation of the average of a matrix element for larger 
tensors.

\proclaim{(4.2) Lemma} Let us fix two tensors 
$A=(a_I)$ and $B=(b_I)$ for \newline
$I=\bigl(1 \leq i_1, \ldots, i_d \leq n\bigr)$. For a positive 
integer $m$ (in particular, for $m=2k$), let us define 
tensors $X=A^{\otimes m}$ and $Y=B^{\otimes m}$ as follows:
$$X=(x_{J}) \quad \text{and} \quad Y=(y_{J}) \quad \text{where} \quad 
J=\bigl(1 \leq j_1, \ldots, j_{dm} \leq n\bigr)$$ 
and where
$$x_J=a_{I_1} \cdots a_{I_{m}} \quad \text{and} \quad 
y_J=b_{I_1} \cdots b_{I_m} \quad \text{provided} \quad 
J=(I_1, \ldots, I_m).$$ 
Then 
$${1 \over n!} \sum_{g \in S_n} 
\big\langle B,\ g A \big\rangle^m =
{1 \over n!} \sum_{g \in S_n} \big\langle Y,\ g X \big\rangle.$$
\endproclaim 
\demo{Proof}
The proof follows by observation that 
$$\big\langle B,\ g A \big\rangle^m =
\big\langle B^{\otimes m},\ g A^{\otimes m} \big\rangle=
\big\langle Y,\ g X \big\rangle.$$
{\hfill \hfill \hfill} \qed
\enddemo
Next, we show how to compute the average.
\proclaim{(4.3) Lemma} Let us fix a positive 
integer $l$ (in particular, $l=md=2kd$). For a partition 
$\Sigma=\bigl\{\Sigma_1, \ldots, \Sigma_r \bigr\}$ of the set 
$\{1, \ldots, l\}$ into non-empty disjoint subsets, we say 
that a sequence $I=(i_1, \ldots, i_l)$ has type $\Sigma$ if 
for each $\Sigma_p$ the indices $i_j: j \in \Sigma_p$ are all equal 
and if for each pair of subsets $\Sigma_p$ and $\Sigma_q$ the 
indices $i_j: j \in \Sigma_p$ and $i_j: j \in \Sigma_q$ are
different. 

Let $X=(x_I)$ and $Y=(y_I)$, 
$I=\bigl(1 \leq i_1, \ldots, i_l \leq n\bigr)$ be tensors
(in particular, we can have $X=A^{\otimes m}=A^{\otimes 2k}$ and 
$Y=B^{\otimes m}=B^{\otimes 2k}$).

Let us define the tensors $\overline{X}=(\overline{x}_I)$ and 
$\overline{Y}=(\overline{y}_I)$,
$I=\bigl(1 \leq i_1, \ldots, i_l \leq n\bigr)$ by 
$$\overline{x}_I={(n-r)! \over n!} 
\sum \Sb J: \text{\ type of\ } J= \text{\ type of\ }
I \endSb 
x_J \quad \text{provided} \quad \text{type\ }I=
\bigl(\Sigma_1, \ldots, \Sigma_r\bigr)$$
and 
$$\overline{y}_I=
{(n-r)! \over n!} \sum \Sb J: \text{\ type \ } J= \text{\ type \ }
I \endSb 
y_J \quad \text{provided} \quad \text{type\ }
I=\bigl(\Sigma_1, \ldots, \Sigma_r\bigr).
$$ 
Then 
$${1 \over n!} \sum_{g \in S_n} \big\langle Y,\ g X \big\rangle
=\big\langle \overline{Y},\ \overline{X} \big\rangle.$$
\endproclaim
\demo{Proof}
The two index sets $I=(i_1, \ldots, i_l)$ and $J=(j_1, \ldots, j_l)$ 
belong to the same orbit of the action $I \longmapsto g I$ of 
$S_n$ if and only if they have the same type
$\bigl\{ \Sigma_1, \ldots, \Sigma_r \bigr\}$. Moreover, the 
stabilizer of $I$ consists of $(n-r)!$ permutations. Hence
$$\overline{X}={1 \over n!} \sum_{g \in S_n} g X 
\quad \text{and} \quad 
\overline{Y}={1\over n!} \sum_{g \in S_n} g Y.$$
We have 
$${1 \over n!} \sum_{g \in S_n} \big\langle Y,\ g X \big\rangle=
\Big\langle \quad {1 \over n!} \sum_{g \in S_n} g Y, 
\quad {1 \over n!} \sum_{g \in S_n} g X \quad \Big\rangle$$
and the proof follows.
{\hfill \hfill \hfill} \qed
\enddemo
Combining Lemmas 4.2 and 4.3, we observe that as long as $d$ and $k$ are 
fixed, we can compute $\|f\|_{2k}$ in $O\bigl(n^{2kd}\bigr)$ time,
that is, in polynomial in $n$ time.

In particular, from Corollary 1.6, we conclude that for any fixed 
$d$ and for any fixed $\epsilon>0$ there exists a polynomial in $n$ algorithm 
for estimating $\|f\|_{\infty}$ within a factor of $\epsilon n^{d/2}$.
This result seems to be new already for $d=2$, cf. \cite{1}.

\remark{(4.4) Remark} So far we have shown how to approximate 
$\|f\|_{\infty}$ by $\|f\|_{2k}$ but we did not discuss how to find
a particular permutation $g$ which gives the value of $|f(g)|$ 
close to $\|f\|_{\infty}$. In fact,
it is not hard to construct a permutation
$g \in S_n$ for which $|f(g)| \geq \|f\|_{2k}$
and hence $|f(g)|$ approximates $\|f\|_{\infty}$ 
within a factor of $\epsilon n^{d/2}$ at the cost of some extra 
work, which still results in a polynomial time algorithm 
when $k$ is fixed. The idea is to use the 
``divide-and-conquer'' approach. We split the symmetric group $S_n$ into the 
union of cosets $S_j=\bigl\{g: g(1)=j \bigr\}$ and then 
compute the average value of $f^{2k}$ over each coset separately 
(this would require some straightforward modification of Lemma 4.3).
Then a coset should be chosen which gives the largest average. 
Thus we have determined $g(1)=j$ and we proceed to 
determine $g(2), \ldots, g(n)$ successively. 
\endremark

\enddocument

\end